\documentclass[letterpaper,10pt,journal]{IEEEtran}

\usepackage{kjkslim}

\usepackage{cite} 

\IEEEoverridecommandlockouts

\title{Convexity and monotonicity in nonlinear optimal control under uncertainty}
\author{Kevin J. Kircher,~\IEEEmembership{Member,~IEEE,} and K. Max Zhang%
\thanks{The authors are with the Sibley School of Mechanical and Aerospace Engineering, Cornell University, Ithaca, NY 14853. {\tt \small \{kjk82,kz33\}@cornell.edu}}}

\begin{document}

\maketitle

\begin{abstract}
We consider the problem of finite-horizon optimal control design under uncertainty for imperfectly observed discrete-time systems with convex costs and constraints. It is known that this problem can be cast as an infinite-dimensional convex program when the dynamics and measurements are linear, uncertainty is additive, and the risks associated with constraint violations and excessive costs are measured in expectation or in the worst case. In this paper, we extend this result to systems with convex or concave dynamics, nonlinear measurements, more general uncertainty structures and other coherent risk measures. In this setting, the optimal control problem can be cast as an infinite-dimensional convex program if (1) the costs, constraints and dynamics satisfy certain monotonicity properties, and (2) the measured outputs can be reversibly `purified' of the influence of the control inputs through Q- or Youla-parameterization. The practical value of this result is that the finite-dimensional subproblems arising in a variety of suboptimal control methods, notably including model predictive control and the Q-design procedure, are also convex for this class of nonlinear systems. Subproblems can therefore be solved to global optimality using convenient modeling software and efficient, reliable solvers. We illustrate these ideas in a numerical example.
\end{abstract}


\section{Introduction}

We consider the problem of operating a system that evolves in discrete time steps over a finite horizon under uncertainty. Our goal is to design an output feedback control policy that minimizes a convex cost while satisfying convex constraints. This problem has been studied in the stochastic setting, where uncertainty is modeled probabilistically, the expected cost is minimized, and constraints are enforced in expectation or in probability. \cite{Whittle1982, Kumar1986, Prekopa1995, Bertsekas2005, Astrom2006, Shapiro2009} It has also been studied in the robust setting, where uncertainty is assumed to come from a given set, the worst-case cost is minimized, and constraints are enforced for all possible realizations of the uncertain influences on the system. \cite{Zhou1996, Ben-Tal1998, Ben-Tal2008, Ben-Tal2009, Bertsimas2011} We accommodate both approaches in this paper, as well as others, adopting the flexible view of risk developed in \cite{Artzner1997, Artzner1999, Rockafellar2000, Rockafellar2002,  Rockafellar2007, Rockafellar2018}.

Optimal control under uncertainty is a hard problem, even when the costs and constraints are convex. This is due in part to the fact that the optimization variable is a \textit{policy}, a collection of functions that map measured system outputs into control inputs at each time step. In the stochastic setting, optimal policies can, in principle, be found analytically using stochastic dynamic programming. \cite{Bellman1957} In practice, analytical solution is typically limited to systems of very low dimension. Notable exceptions are linear systems with additive uncertainty, no constraints, and costs that are either quadratic \cite{Newton1957, Kalman1960} or exponential-of-quadratic \cite{Jacobson1973, Whittle1981}.

When analytical solution is impractical, various methods can generate suboptimal solutions that often perform well. Some examples are classical linear feedback control design, reinforcement learning \cite{Bertsekas1995, Sutton1998, Lillicrap2015, Mnih2016} and approximate dynamic programming \cite{Bertsekas2005, DeFarias2003, Powell2007, Lewis2013}, approximation methods for multistage stochastic programming \cite{Shapiro2009, Swamy2005, Gupta2005, Donohue2006, Kuhn2011} and robust optimization with recourse \cite{Kuhn2011, Ben-Tal2006, Beyer2007, ChenZhang2009}, the $Q$-design procedure \cite{Kucera1980, Desoer1984, Anantharam1984, Boyd1991, Skaf2009, Skaf2010}, and model predictive control (MPC) in its certainty-equivalent \cite{Mayne2000, Bemporad2002, Wang2010}, stochastic \cite{Schwarm1999, Oldewurtel2008, Blackmore2010, Cannon2011, Cinquemani2011, Schildbach2014, KircherSA2016}, and robust \cite{Campo1987, Magni2003, Mayne2005, Goulart2007, CalafioreFagiano2013} variants. When perfect state information is not available, control methods may be paired with a state estimator such as a linear \cite{KalmanKF1960, Kalman1961}, extended \cite{Kopp1963, Cox1964, Ljung1979} or unscented \cite{Wan2000, VanDerMerwe2001} Kalman filter or a particle filter \cite{Kitagawa1996, Arulampalam2002}.

These methods vary widely in their scope, scalability and performance. A common theme, however, is that they tend to work best for linear systems. This is due in part to the fact that for linear systems, (an equivalent transformation of) the optimal control problem is convex. \cite{Skaf2009} Suboptimal control methods often involve numerically solving optimization subproblems generated by the original optimal control problem. The convexity of the original problem typically carries over to the subproblems, allowing them to be efficiently and reliably solved to global optimality. When subproblems are nonconvex, however, guarantees of global optimality are generally unavailable, and solvers and initial guesses may need to be carefully tailored to the applications at hand.

In \cite{Skaf2009}, Skaf and Boyd demonstrate that for linear systems with additive uncertainty, the control design problem can be transformed to an infinite-dimensional convex program. Their method hinges on a change of variables related to the $Q$- or Youla-parameterization \cite{Youla1976-1, Youla1976-2, Wu2010} and to purified output feedback control \cite{Kumar1986, Ben-Tal2009, Ben-Tal2006}. When perfect state information is available, this change of variables parameterizes state feedback policies by equivalent disturbance feedback policies. Similar arguments have justified the use of (typically affine) disturbance feedback policies in robust and stochastic MPC of perfectly-observed linear systems. \cite{Skaf2010, Oldewurtel2008, VanHessem2002, Lofberg2003, Ben-Tal2004, Goulart2006, WangOng2010}

For nonlinear systems, the optimal control problem is widely understood to be nonconvex due to nonlinear equality constraints introduced by nonlinear dynamics. Equality constraints can be eliminated, however, by iteratively applying the dynamics to express the state trajectory in terms of the control and exogenous input trajectories. In \cite{Rantzer2014}, Rantzer and Bernhardsson observe that in convex-monotone systems, where the dynamics are convex and nondecreasing in the state and control input, the state trajectory is \textit{convex} in the control input trajectory. In \cite{Schmitt2017}, Schmitt \etal \ generalize this result to convex-state-monotone systems, where the dynamics need not be monotone in the control input. 

An immediate consequence of the observations in \cite{Rantzer2014, Schmitt2017} is that for some nonlinear systems, \textit{open-loop} optimal control (where the decision variable is a fixed control trajectory, rather than a feedback policy) is a convex optimization problem. This holds for convex-state-monotone systems in particular, provided the cost and constraints are nondecreasing in the states. This raises two further questions: Are there other nonlinear systems for which open-loop optimal control is convex? What about \textit{closed-loop} optimal control, where we optimize over policies?

This discussion motivates the definition of a \textit{convex system} as one for which open-loop optimal control is a convex optimization problem. After setting the stage in \S\ref{problem}, we establish three results for convex systems in this paper.
\ben
\item \textbf{Characterization (\S\ref{convex-systems}):} Systems with mixed convex and linear dynamics are convex systems, provided (a) any linear dynamics are independent of states with nonlinear dynamics, and (b) the cost, constraints and nonlinear dynamics are nondecreasing in the states with nonlinear dynamics. Concave dynamics can also be accommodated.
\item \textbf{Convex closed-loop design (\S\ref{convexity}):} If the measured outputs can be reversibly `purified' of the influence of the control inputs, then the closed-loop optimal control problem for convex systems can be transformed to an equivalent infinite-dimensional convex program. The transformation involves changing variables to policies in the purified outputs via $Q$- or Youla-parameterization.
\item \textbf{Approximate solution (\S\ref{approximation}--\ref{example}):} The finite-dimensional subproblems arising in a variety of suboptimal control methods, notably including MPC and the $Q$-design procedure, are convex for convex systems. Subproblems can therefore be solved to global optimality using convenient modeling software and efficient, reliable solvers. We illustrate this in a numerical example.
\een

\section{Problem statement}
\label{problem}

\subsection{System}

We consider a system to be operated over a finite discrete time span $t = 0, \dots, T$. We can influence the system through the control inputs $u_0, \dots, u_{T-1}$. The system is also influenced by exogenous inputs $\delta_0, \dots, \delta_{T}$, which are generally uncertain.\footnote{We view the exogenous input $\delta_t$ as a random vector defined on an underlying probability space. We do not assume that the joint distribution of the exogenous input trajectory is known, but we do assume this distribution is independent of the control input trajectory.} The exogenous inputs may include process disturbances, sensor noise, initial states, uncertain model parameters, prices, command or reference signals, \etc \ The exogenous inputs could come from bounded or unbounded sets, and need not be independent or identically distributed over time.

The control and exogenous inputs determine the states $x_0, \dots, x_T$ through the system dynamics:
\bneq
\begin{aligned}
x_0 &= f_0(\delta_0) \\
x_t &= f_t(x_{t-1}, u_{t-1}, \delta_t), \quad t = 1, \dots, T . \label{dynamics}
\end{aligned}
\eneq
We observe the system through the measured outputs
\bneq
\begin{aligned}
y_0 &= g_0(x_0,\delta_0) \\
y_t &= g_t(x_t,u_{t-1},\delta_t), \quad t = 1, \dots, T-1 . \label{measurements}
\end{aligned}
\eneq
We assume that the dynamics mappings $f_0, \dots, f_T$ and measurement mappings $g_0, \dots, g_{T-1}$ are known. At each time $t$, the controller receives the measured output $y_t$. It decides the control input
\[
u_t = \pi_t(y_0, \dots, y_t) 
\]
by evaluating an output feedback control law $\pi_t$. The control policy $\pi = (\pi_0, \dots, \pi_{T-1})$ is designed in advance; this design problem is the subject of this paper.

%

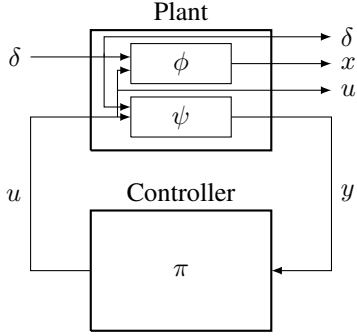
\begin{figure}
\centering
\begin{tikzpicture}[scale=0.8]
	\pgfmathsetmacro{\wc}{3}
	\pgfmathsetmacro{\hc}{2}
	\pgfmathsetmacro{\yp}{3}
	\pgfmathsetmacro{\dy}{0.222}
	\pgfmathsetmacro{\dx}{0.666}
	\pgfmathsetmacro{\hm}{0.666}
	\pgfmathsetmacro{\wsig}{1}
	
	\draw[thick] (0,0) -- (\wc,0) -- (\wc,\hc) -- (0,\hc) -- (0,0);
	\node at (\wc/2,\hc/2) {$\pi$};
	\node[above] at (\wc/2,\hc) {Controller};
	
	\draw[thick] (0,\yp) -- (\wc,\yp) -- (\wc,\yp+\hc) -- (0,\yp+\hc) -- (0,\yp);
	\node[above] at (\wc/2,\yp+\hc) {Plant};

	\draw (\dx,\yp+\dy) -- (\wc-\dx,\yp+\dy) -- (\wc-\dx,\yp+\dy+\hm) -- (\dx,\yp+\dy+\hm) -- (\dx,\yp+\dy);
	\node at (\wc/2,\yp+\dy+\hm/2) {$\psi$};
	
	\draw (\dx,\yp+2*\dy+\hm) -- (\wc-\dx,\yp+2*\dy+\hm) -- (\wc-\dx,\yp+2*\dy+2*\hm) -- (\dx,\yp+2*\dy+2*\hm) -- (\dx,\yp+2*\dy+\hm);
	\node at (\wc/2,\yp+2*\dy+1.5*\hm) {$\phi$};
	
	\draw[->] (\wc-\dx,\yp+\dy+\hm/2) -- (\wc+\wsig,\yp+\dy+\hm/2) -- (\wc+\wsig,\hc/2) -- (\wc,\hc/2);
	\node[above right] at (\wc+\wsig,\hc) {$y$};

	\draw[->]  (0,\hc/2) -- (-\wsig,\hc/2) -- (-\wsig,\yp+\dy+\hm/2) -- (\dx,\yp+\dy+\hm/2);
	\node[above left] at (-\wsig,\hc) {$u$};

	\draw[->] (-\wsig,\yp+2*\dy+\hm+2*\hm/3) -- (\dx,\yp+2*\dy+\hm+2*\hm/3);
	\node[left] at (-\wsig,\yp+2*\dy+\hm+2*\hm/3) {$\delta$};
	
	\draw[->] (\wc-\dx,\yp+2*\dy+3*\hm/2) -- (\wc+\wsig,\yp+2*\dy+3*\hm/2);
	\node[right] at (\wc+\wsig,\yp+2*\dy+3*\hm/2) {$x$};
	
	\draw[->] (2*\dx/3,\yp+\dy+\hm/2) -- (2*\dx/3,\yp+2*\dy+\hm+\hm/3) -- (\dx,\yp+2*\dy+\hm+\hm/3);
	
	\draw[->] (2*\dx/3,\yp+\hc/2) -- (\wc+\wsig,\yp+\hc/2);
	\node[right] at (\wc+\wsig,\yp+\hc/2) {$u$};
	
	\draw[->] (\dx/3,\yp+2*\dy+\hm+3*\hm/4) -- (\dx/3,\yp+\dy+3*\hm/4) -- (\dx,\yp+\dy+3*\hm/4);
	
	\draw[->] (\dx/3,\yp+2*\dy+\hm+2*\hm/3) -- (\dx/3,\yp+\hc-\dy/2) -- (\wc+\wsig,\yp+\hc-\dy/2);
	\node[right] at (\wc+\wsig,\yp+\hc-\dy/2) {$\delta$};
	
\end{tikzpicture}
\caption{The closed-loop system is driven by the exogenous input $\delta$. The control input $u$ is determined by the measured output $y$ through the control policy $\pi$. Performance is evaluated based on $\delta$, $u$ and the state $x$.}
\label{systemFig}
\end{figure}

To simplify specification of the policy design problem, we work with the input, state and output trajectories. By iteratively applying the dynamics, the state trajectory can be expressed in terms of the control and exogenous input trajectories as
\[
x = \bmat
x_0 \\
x_1 \\
\vdots \\
x_T \\
\emat = \bmat 
\phi_0(\delta_{0}) \\
\phi_1(u_0,\delta_{0:1}) \\
\vdots \\
\phi_T(u_{0:T-1},\delta_{0:T}) \\
\emat = \phi(u,\delta) .
\]
Here the subscript $t_1:t_2$ denotes a trajectory from time $t_1$ to $t_2$. For example, $\delta_{0:T} = (\delta_0, \dots, \delta_T) = \delta$ and $u_{0:T-1} = (u_0,\dots,u_{T-1}) = u$. 
The input-state mappings $\phi_t$ are defined by the recursion
\[
\begin{aligned}
\phi_0(\delta_{0}) &= f_0(\delta_0) \\
\phi_1(u_0, \delta_{0:1}) &= f_1(\phi_0(\delta_0),u_0,\delta_1) \\
\phi_t(u_{0:t-1},\delta_{0:t}) &= f_t(\phi_{t-1}(u_{0:t-2},\delta_{0:t-1}),u_{t-1},\delta_t), \\
&\qquad t = 2, \dots, T.
\end{aligned}
\]
Similarly, the measured output trajectory can be written as
\[
y = \bmat
y_0 \\
y_1 \\
\vdots \\
y_{T-1} \\
\emat = \bmat 
\psi_0(\delta_{0}) \\
\psi_1(u_0,\delta_{0:1}) \\
\vdots \\
\psi_{T-1}(u_{0:T-2},\delta_{0:T-1}) \\
\emat = \psi(u,\delta) ,
\]
where the input-output mappings $\psi_t$ are defined recursively by
\[
\begin{aligned}
\psi_0(\delta_0) &= g_0(f_0(\delta_0),\delta_0) \\
\psi_t(u_{0:t-1},\delta_{0:t}) &= g_t(\phi_t(u_{0:t-1},\delta_{0:t}), u_{t-1}, \delta_{t}), \\
&\qquad t = 1, \dots, T-1 .
\end{aligned}
\]
We also write the control input trajectory as
\[
u = \pi(y) = \bmat
\pi_0(y_0) \\
\pi_1(y_{0:1}) \\
\vdots \\
\pi_{T-1}(y_{0:T-1}) \\
\emat .
\]
Figure \ref{systemFig} illustrates the system.

\subsection{Cost, constraints and risk measures}
\label{risk}

We are interested in designing the output feedback control policy $\pi = (\pi_0, \dots, \pi_{T-1})$. The policy must be \textit{causal}, meaning each control law can depend on outputs from the past or present, but not the future. We denote the set of causal output feedback policies by $\Pi$.  Our goal is to find a $\pi \in \Pi$ that minimizes a cost
\bneq
c_0(x,u,\delta) \label{cost}
\eneq
while satisfying constraints
\bneq
c_j(x,u,\delta) \leq 0, \quad j = 1, \dots, J. \label{constraints}
\eneq
We assume that each scalar-valued function $c_0, \dots, c_J$ is convex in $(x,u)$ for all $\delta$. 

Because the exogenous inputs are uncertain, the goals of minimizing cost and satisfying constraints are ambiguous. Should they be accomplished on average with respect to the distribution of $\delta$, for all possible realizations of $\delta$, or somewhere in between? Resolving these ambiguities amounts to choosing measures of risk.\footnote{For more on measuring risk, a rich subject that has received much recent attention, we refer the reader to \cite{Artzner1997, Artzner1999, Rockafellar2000, Rockafellar2002,  Rockafellar2007, Rockafellar2018}.} A \textit{risk measure} $\risk$ is a functional that maps uncertain scalars into deterministic scalars (or possibly $\infty$). It quantifies the risk associated with positive realizations of its uncertain argument. In the language of \cite{Artzner1997, Artzner1999}, the risk measure $\risk_j$ should be associated with constraint $j$ if the designer views the risk of constraint violation as \textit{acceptable} whenever
\[
\risk_j c_j(x,u,\delta) \leq 0 .
\]
The risks associated with excessive costs can be treated similarly to the risks associated with constraint violations, because minimizing $c_0(x,u,\delta)$ is equivalent to minimizing a deterministic scalar $\tau$ subject to the constraint $c_0(x,u,\delta) - \tau \leq 0$.

\begin{table}
\centering
\caption{Example risk measures for a random variable $X$.}
\smallskip
\begin{tabular}{l | l | l }
Risk measure 			&Definition of $\risk X$			&Risk acceptable if\dots \\
\hline
Expected value			&$\Expect X$					& $X \leq 0$ on average \\
Predicted value			&prediction $\hat X$ of $X$		& $\hat X \leq 0$ \\
Worst-case value		&$\sup X$						& $X \leq 0$ almost surely \\
\begin{tabular}[x]{@{}l@{}} Value-at-risk,\\$\quad \VaR_\beta(X)$ \end{tabular}	&\begin{tabular}[x]{@{}l@{}} $\inf \{\alpha \mid $\\$\quad \Prob \setof{X \leq \alpha} \geq \beta \}$ \end{tabular}	& $\Prob \setof{X \leq 0} \geq \beta$ \\
\begin{tabular}[x]{@{}l@{}} Conditional \\\quad value-at-risk \end{tabular}	&\begin{tabular}[x]{@{}l@{}} $\inf_\alpha \{\alpha \ +$\\$\quad \Expect (X - \alpha)_+ / (1-\beta) \}$ \end{tabular} & \begin{tabular}[x]{@{}l@{}} $\Expect[X \mid X \geq $\\$\quad \VaR_\beta(X)] \leq 0$ \end{tabular} \\
\end{tabular}
\label{risk-table}
\end{table}

Table \ref{risk-table} describes five risk measures commonly used in optimization under uncertainty. The expected value is risk-neutral: it weighs constraint violations and slacks equally, modeling the risk of constraint violation as acceptable if slacks balance out violations on average. The worst-case value, by contrast, is maximally risk-averse: it models any possibility of constraint violation as unacceptable. The value-at-risk at confidence level $\beta \in [0,1]$ ($\VaR_\beta$) is the smallest threshold that, with 100$\beta$\% confidence, will not be exceeded. Typical values of $\beta$ are 0.9, 0.95 and 0.99. $\VaR_\beta$ generates chance constraints of the form
\[
\Prob \setof{c_j(x,u,\delta) \leq 0} \geq \beta .
\]
The conditional value-at-risk at confidence level $\beta$ ($\cvar_\beta$) is the conditional expectation, conditioned on the event that $\VaR_\beta$ is exceeded. \cite{Rockafellar2000, Rockafellar2002} $\cvar_\beta$ always upper bounds $\VaR_\beta$. The $\cvar_\beta$ of $c_j(x,u,\delta)$ can be constrained by minimizing a scalar $\alpha$ subject to
\[
\Expect (c_j(x,u,\delta) - \alpha)_+ \leq - \alpha (1-\beta) ,
\]
where $(\cdot)_+$ denotes the positive part, $\max \setof{0, \cdot}$. An optimal $\alpha^\star$ upper bounds $\VaR_\beta c_j(x,u,\delta)$; this bound is usually tight. \cite{Rockafellar2002}

In this paper, we leave the choice of risk measures to the designer. This flexibility makes the results in the following sections applicable in a range of settings, including the robust and stochastic frameworks and others. We do, however, impose one restriction on the risk measures: we assume that they are \textit{coherent} in the sense of \cite{Artzner1997, Artzner1999, Rockafellar2007, Rockafellar2018}. By this, we mean that each risk measure is convex, nondecreasing, lower semicontinuous, and preserves certainty (\ie, $\risk (C) = C$ for constant $C$). We make this assumption primarily because coherent risk measures preserve convexity: if $\risk_j$ is coherent and $(x,u) \mapsto c_j(x,u,\delta)$ is convex for all $\delta$, then $(x,u) \mapsto \risk_j c_j(x,u,\delta)$ is convex. \cite{Rockafellar2007} 

All but one of the risk measures mentioned above are coherent. The exception is $\VaR_\beta$, which is convex only in structured special cases. We note, however, that $\VaR_\beta$ can be approximated in a convexity-preserving fashion using sampling \cite{Calafiore2005, Campi2008, Calafiore2010, Schildbach2013, Zhang2015, Campi2018} or conservative upper bounds such as $\cvar_\beta$ \cite{Shapiro2009, Cinquemani2011}.

In summary, we are interested in the following problem of optimal control under uncertainty:
\bneq
\begin{array}{ll}
\underset{\pi \in \Pi}{\minimize}	& \risk_0 c_0(x, u, \delta)  \\ \label{OC-pi} \tag{OC}
\mbox{subject to}	& \risk_j c_j(x, u, \delta) \leq 0 , \quad j = 1, \dots, J \\
				& x = \phi(u, \delta) \\
				& u = \pi(y) \\
				& y = \psi(u, \delta) .
\end{array}
\eneq
This problem is intractable in its current form. This is due in part to the fact that the decision variable $\pi$ is infinite-dimensional, an obstacle that we will address approximately in \S\ref{approximation}. Problem \eqref{OC-pi} is also complicated by the interdependence between the state, output and control trajectories. We will unravel this interdependence in \S\ref{convexity} using a nonlinear change of variables. First, we will consider the case of open-loop control, where the control laws $\pi_t$ are restricted to be constant.

\section{Convex systems}
\label{convex-systems}

To understand when closed-loop control design is a convex optimization problem, we begin by building intuition in the simpler context of open-loop control. In this context, a (finite-dimensional) vector of control inputs $\tilde u = (\tilde u_0, \dots, \tilde u_{T-1})$  is decided in advance and implemented without feedback. This amounts to drastically restricting the search space from the set of all causal output feedback policies to the subspace of \textit{constant} policies. With this restriction, the open-loop optimal control problem is to
\bneq
\begin{array}{ll}
\underset{\tilde u}{\minimize}	& \risk_0 c_0(\phi(\tilde u, \delta), \tilde u,  \delta)  \\ \label{OC-OL} \tag{OLOC}
\mbox{subject to}			& \risk_j c_j(\phi(\tilde u, \delta), \tilde u,  \delta) \leq 0 , \quad j = 1, \dots, J .
\end{array}
\eneq
This is a finite-dimensional optimization problem in the vector $\tilde u$. Because coherent risk measures preserve convexity \cite{Rockafellar2007}, Problem \eqref{OC-OL} is convex if each cost and constraint function $c_j$, when composed with the input-state mapping $\phi$, is convex in $u$ for all $\delta$. This motivates the following definition.

\begin{definition}[Convex system] We call the system with dynamics \eqref{dynamics}, cost \eqref{cost} and constraints \eqref{constraints} a {\rm convex system} if the functions
\[
u \mapsto c_j(\phi(u, \delta),u, \delta), \quad j = 0, \dots, J
\]
are convex for all $\delta$. 
\end{definition}
\noindent By definition, open-loop optimal control of convex systems is a convex optimization problem. In \S\ref{convexity}, we will show that \textit{closed-loop} optimal control of convex systems is also convex, provided the outputs can be reversibly `purified' of the influence of the control inputs. First, we will characterize a class of convex systems.

\subsection{Characterizing convex systems}

To understand when a system is convex, we need to understand how the states and controls propagate through the dynamics \eqref{dynamics} into the cost \eqref{cost} and constraints \eqref{constraints}. If the dynamics are linear, then this process is straightforward. The state trajectory $x$ is affine in the control trajectory $u$ for all $\delta$. Convexity is preserved under composition with an affine mapping \cite{Boyd2004}, so for linear systems $u \mapsto c_j(\phi(u, \delta),u, \delta)$ is convex for all $\delta$. 

When the system is nonlinear, the process is less straightforward. We do not explore it in general here. Instead, we consider a class of systems with the following dynamics:
\bneq
\begin{aligned}
x_t = \bmat
x^{\rm aff}_{t} \\
x^{\rm cvx}_{t} \\
\emat &= \bmat
A_t(\delta_t) x^{\rm aff}_{t-1} + B_t(\delta_t) u_{t-1} + w_t(\delta_t) \\
h_{t}(x^{\rm aff}_{t-1}, x^{\rm cvx}_{t-1}, u_{t-1},\delta_{t}) \\ 
\emat . \label{linear-convex}
\end{aligned}
\eneq
For this class of systems, convexity can be established using simple composition rules. An important restriction is that the states $x^\text{\rm aff}_t$ with linear dynamics are independent of the nonlinear states $x^{\rm cvx}_t$. This ensures that the trajectory $x^\text{\rm aff} = (x^\text{\rm aff}_0, \dots, x^\text{\rm aff}_T)$ is affine in the control trajectory. Any states with linear dynamics that depend on nonlinear states are included in $x^{\rm cvx}_t$.

\begin{theorem} \label{convex-system-theorem} The system with dynamics \eqref{linear-convex}, cost \eqref{cost} and constraints \eqref{constraints} is a convex system if the following conditions hold for all $\delta$.
\ben
\item For $t = 1, \dots, T$, each row of the nonlinear dynamics mapping $h_t$ is 
	\ben
	\item jointly convex in $(x^{\rm aff}_{t-1}, x^{\rm cvx}_{t-1}, u_{t-1})$, and
	\item nondecreasing in each element of $x^{\rm cvx}_{t-1}$.
	\een
\item For $j = 0, \dots, J$, the function
\[
(x^{\rm aff},x^{\rm cvx},u) \mapsto c_j((x^{\rm aff},x^{\rm cvx}),u,\delta)
\]
is nondecreasing in each element of $x^{\rm cvx} = (x^\text{\rm cvx}_0, \dots, x^\text{\rm cvx}_T)$.
\een
\end{theorem}
\noindent Appendix \ref{convex-system-proof} contains a proof of Theorem \ref{convex-system-theorem}. The proof hinges on the fact that the composition of a convex nondecreasing function with a convex function is convex. \cite{Boyd2004}

Three remarks are in order. First, system \eqref{linear-convex} is a straightforward generalization of the convex-state-monotone system studied in \cite{Schmitt2017}, which in turn generalizes the convex-monotone system studied in \cite{Rantzer2014}. A precise name for system \eqref{linear-convex} would be \textit{convex-nonlinear-state-monotone}, since the monotonicity requirement applies only to the states with nonlinear dynamics. The purpose of this generalization is to absorb the class of linear systems for which optimal control was shown in \cite{Skaf2009} to be convex. In particular, the linear-quadratic system is not convex-monotone or convex-state-monotone in general, but is a special case of system \eqref{linear-convex} with linear dynamics, additive uncertainty, no constraints and quadratic cost.

Second, concave dynamics can be included in system \eqref{linear-convex} after a sign change. For example, we consider a system with
\[
\begin{aligned}
z^\text{cvx}_{t} &= h^\text{cvx}_{t}(x^{\rm aff}_{t-1}, z^\text{cvx}_{t-1}, z^\text{ccv}_{t-1}, u_{t-1}, \delta_t) \\
z^\text{ccv}_{t} &= h^\text{ccv}_{t}(x^{\rm aff}_{t-1}, z^\text{cvx}_{t-1}, z^\text{ccv}_{t-1}, u_{t-1}, \delta_t) ,
\end{aligned}
\]
where each row of $h^\text{cvx}_{t}$ is nondecreasing in $z^\text{cvx}_{t-1}$, nonincreasing in $z^\text{ccv}_{t-1}$ and convex. Similarly, each row of $h^\text{ccv}_{t}$ is nondecreasing in $z^\text{ccv}_{t-1}$, nonincreasing in $z^\text{cvx}_{t-1}$ and concave. If we define the nonlinear state as
\[
\begin{aligned}
x^{\rm cvx}_t = \bmat
z^\text{cvx}_{t} \\
- z^\text{ccv}_{t} \\
\emat &= \bmat
h^\text{cvx}_{t}(x^{\rm aff}_{t-1}, z^\text{cvx}_{t-1}, z^\text{ccv}_{t-1}, u_{t-1}, \delta_t) \\
- h^\text{ccv}_{t}(x^{\rm aff}_{t-1}, z^\text{cvx}_{t-1}, z^\text{ccv}_{t-1}, u_{t-1}, \delta_t) \\
\emat \\
&= h_t(x^{\rm aff}_{t-1}, x^{\rm cvx}_{t-1}, u_{t-1}, \delta_t) ,
\end{aligned} 
\]
then each row of $h_t$ is nondecreasing in $x^{\rm cvx}_{t-1}$ and convex, as required by Theorem \ref{convex-system-theorem}. 

Third, system \eqref{linear-convex} includes a rich class of nonlinear systems. Provided careful attention is paid to curvature and monotonicity, it can admit exponentials and logarithms, quadratic forms, roots, powers, nonsmooth functions such as maxima, minima and absolute values, and sums and compositions of the above. Many more examples of nonlinear convex and concave functions can be found in \S3 of \cite{Boyd2004}. Systems with nonlinear convex dynamics have arisen in applications ranging from cancer and HIV treatment scheduling \cite{Rantzer2014} to voltage control in power systems \cite{Rantzer2014} to freeway congestion management \cite{Schmitt2017} to energy storage control \cite{Khodabakhsh2016, Xu2017}. The drug treatment model in \cite{Rantzer2014} is a special case of a more general class of bilinear systems that, through a logarithmic transformation, can be cast as convex systems.

\section{Convex closed-loop optimal control}
\label{convexity}

By contrast to the open-loop optimal control problem \eqref{OC-OL}, the equality constraints in the closed-loop optimal control problem \eqref{OC-pi} cannot be easily eliminated. This is due to the complicated interdependence between the state, output and control input trajectories. Under some conditions on the uncertainty structure, however, this interdependence can be disentangled through a nonlinear change of variables related to the $Q$- or Youla-parameterization \cite{Youla1976-1, Youla1976-2, Wu2010} and purified output feedback control \cite{Kumar1986, Ben-Tal2009, Ben-Tal2006}. We now establish sufficient conditions for this change of variables to be possible.

\subsection{Purifiability and Q-parameterization}
\label{purifiable-Q}

\begin{definition}[Purifiable] We call the system with dynamics \eqref{dynamics} and measurements \eqref{measurements} {\rm purifiable} if for each $t = 1, \dots, T-1$, there exist mappings $p_t$, $q_t$ and $\xi_t$  such that for all $u_{0:t-1}$ and $\delta_{0:t}$,
\[
\begin{aligned}
&p_t(\psi_{0:t}(u_{0:t-1},\delta_{0:t}), u_{0:t-1}) = \xi_t(\delta_{0:t}) \\
\iff \quad &q_t(\xi_{0:t}(\delta_{0:t}), u_{0:t-1}) = \psi_{t}(u_{0:t-1},\delta_{0:t}) .
\end{aligned}
\]
We call $e_t = \xi_t(\delta_{0:t})$ the purified output and $p_t$ the purifier. 
\end{definition}
\noindent We note for future reference that $y_0$ is a function of $\delta_0$ only:
\[
y_0 = g_0(x_0,\delta_0) = g_0(f_0(\delta_0),\delta_0)  .
\]
Without loss of generality, therefore, we define $p_0$, $q_0$ and $\xi_0$ such that
\[
p_0(y_0) = y_0, \quad q_0(e_0) = e_0, \quad \xi_0(\delta_0) = g_0(f_0(\delta_0),\delta_0) .
\]

\begin{theorem} \label{purifiability-Q-theorem} If the system with dynamics \eqref{dynamics} and measurements \eqref{measurements} is purifiable, then there exists a one-to-one correspondence between causal output feedback policies $\pi$ and causal policies $Q$ in the purified output $e = \xi(\delta)$. Furthermore, given a causal $Q$ in $e$, the unique corresponding causal $\pi$ in $y$ can be constructed from the following recursion:
\bneq
\begin{aligned}
\pi_0(y_0) &= Q_0(y_0) \\ 
\pi_t(y_{0:t}) &= Q_t(p_{0:t}(y_{0:t}, \pi_{0:t-1}(y_{0:t-1}))), \\
&\qquad t = 1, \dots, T-1 . \label{Q-recursion}
\end{aligned}
\eneq
\end{theorem}
\noindent Theorem \ref{purifiability-Q-theorem} is closely related to the nonlinear discrete-time Youla parameterization presented by Wu and Lall in \cite{Wu2010}. We include a proof in Appendix \ref{purifiability-Q-proof} for completeness. We give some examples of purifiable nonlinear systems in \S\ref{purifiable-examples}.

Theorem \ref{purifiability-Q-theorem} establishes that if the system is purifiable, then we can optimize over policies $Q$ in the purified output $e = \xi(\delta)$. We interpret $e_t$ as what remains of the output $y_t$ when the effect of the control inputs $u_0, \dots, u_{t-1}$ has been removed. Given a causal $Q$, the unique corresponding causal $\pi$ can be recovered. We denote the set of causal $Q$ by $\mathcal Q$.

Under this change of variables, an equivalent reformulation of the closed-loop optimal control problem \eqref{OC-pi} is to
\[
\begin{array}{ll}
\underset{Q \in \mathcal Q}{\minimize}	& \risk_0 c_0(x, u, \delta)  \\
\mbox{subject to}	& \risk_j c_j(x, u, \delta) \leq 0 , \quad j = 1, \dots, J \\
				& x = \phi(u, \delta) \\
				& u = Q(e) \\
				& e = \xi(\delta) .
\end{array}
\]
The equality constraints can now be eliminated, giving another equivalent problem:
\bneq
\begin{array}{ll}
\underset{Q \in \mathcal Q}{\minimize}	& \risk_0 c_0(\phi(Q(\xi(\delta)), \delta), Q(\xi(\delta)), \delta)  \\ \label{OC-Q} \tag{OC-$Q$}
\mbox{subject to}	& \risk_j c_j(\phi(Q(\xi(\delta)), \delta), Q(\xi(\delta)), \delta) \leq 0 , \\
&\qquad j = 1, \dots, J .
\end{array}
\eneq
Problem \eqref{OC-Q} is structurally identical to the open-loop optimal control problem \eqref{OC-OL}, except that the optimization is over (infinite-dimensional) purified output feedback policies $Q$ rather than (finite-dimensional) control input trajectories $\tilde u$. It follows that the two problems are convex under the same conditions, namely for convex systems. This observation, together with Theorem \ref{purifiability-Q-theorem}, gives the following result.

\begin{corollary} \label{convex-corollary} For purifiable systems, the output feedback policy design problem \eqref{OC-pi} and the purified output feedback policy design problem \eqref{OC-Q} are equivalent. For convex systems, Problem \eqref{OC-Q} is convex.
\end{corollary}

We will see in \S\ref{approximation} that Corollary \ref{convex-corollary} establishes the basic tractability of a variety of suboptimal control schemes for nonlinear convex systems. First, we provide a few examples of purifiable systems.

\subsection{Purifiability examples}
\label{purifiable-examples}

Purifiability is essentially an invertibility property on the input-output mapping $\psi$ with respect to the exogenous inputs. The requirements for purifiability are (1) at each time step, the influence of the control input history can be removed from the current output, possibly using the output history; and (2) this process must be reversible, in the sense that the current output can be reconstructed from the purified output history and the control input history. To ground this notion, we now provide some concrete examples of purifiable systems. This list is not exhaustive.

\paragraph{Measured exogenous inputs} If the exogenous inputs are measured exactly ($y_t = \delta_t$), then the system is purifiable with
\[
\begin{aligned}
e_t &= \delta_t \\
p_t(y_{0:t}, u_{0:t-1}) &= y_t \\
q_t(e_{0:t},u_{0:t-1}) &= e_t .
\end{aligned}
\]

\paragraph{Pure estimation} If the controller can observe the system but not influence it, then the states and outputs can be expressed as
\[
\begin{aligned}
x_t &= f_t(x_{t-1}, \delta_t) = \phi_t(\delta_{0:t}) \\
y_t &= g_t(x_t, \delta_t) = g_t(\phi_t(\delta_{0:t}), \delta_t) = \psi_t(\delta_{0:t}) .
\end{aligned}
\]
In this case, the system is trivially purifiable with $e_t = y_t$. This implies that various constrained estimation problems can be put in the form of Problem \eqref{OC-Q}. For example, we consider the problem of designing a state estimator $\pi$ such that $u_t = \pi_t(y_{0:t})$ minimizes the mean squared error in estimating $x_t$, with the prior knowledge that $x_t \succeq 0$ almost surely. This can be put in the form of Problem \eqref{OC-Q} by setting $c_0(x,u) = \norm{x-u}_2^2$, $\risk_0 = \Expect$, $c_1(u) = \max_{i,t} \setof{ -(u_t)_i }$ and $\risk_1 = \sup$.

\paragraph{Perfect state information, invertible dynamics} If the states are measured exactly ($y_t = x_t$) and the dynamics are invertible in the exogenous inputs, \ie, there exist mappings $f_t^{-1}$ such that
\[
f_t^{-1}(f_t(x_{t-1},u_{t-1},\delta_t),x_{t-1},u_{t-1}) = \delta_t ,
\]
then the system is purifiable with
\[
\begin{aligned}
e_t &= \delta_t \\
p_t(y_{0:t},u_{0:t-1}) &= f_t^{-1}(y_t,y_{t-1},u_{t-1}) \\
q_t(e_{0:t},u_{0:t-1}) &= \phi_t(u_{0:t-1}, e_{0:t}) .
\end{aligned}
\]
In the special case of additive disturbances, $x_t = f_t(x_{t-1}, u_{t-1}) + \delta_t$, the purifier reduces to 
\[
p_t(y_{0:t},u_{0:t-1}) = y_t - f_t(y_{t-1},u_{t-1}) .
\]

\paragraph{Deterministic dynamics, invertible measurements} If the initial state is measured exactly ($y_0 = x_0$), the dynamics are deterministic,
\[
x_t = f_t(x_{t-1},u_{t-1}) = \phi_t(u_{0:t-1}, x_0) ,
\]
and the measurements are invertible in the exogenous inputs, \ie, there exist mappings $g_t^{-1}$ such that
\[
g_t^{-1}(g_t(x_{t},u_{t-1},\delta_t), x_t, u_{t-1}) = \delta_t ,
\]
then the system is purifiable with
\[
\begin{aligned}
e_t &= \delta_t \\
p_t(y_{0:t}, u_{0:t-1}) &= g_t^{-1}(y_t, \phi_t(u_{0:t-1},y_0), u_{t-1}) \\
q_t(e_{0:t},u_{0:t-1}) &= g_t(\phi_t(u_{0:t-1}, e_{0}), u_{t-1}, e_{t}) .
\end{aligned}
\]
In the special case of additive noise, $y_t = g_t(x_t,u_{t-1}) + \delta_t$, the purifier reduces to 
\[
p_t(y_{0:t}, u_{0:t-1}) = y_t - g_t(\phi_t(u_{0:t-1},y_0),u_{t-1}) .
\]

\paragraph{State-affine dynamics and measurements, additive uncertainty} If the dynamics and measurements have the form
\[
\begin{aligned}
f_t(x_{t-1},u_{t-1},\delta_t) &= A_t x_{t-1} + f_t^u( u_{t-1} ) + w_t(\delta_t) \\
g_t(x_{t},u_{t-1},\delta_t) &= C_t x_{t} + g_t^u( u_{t-1} ) + v_t(\delta_t) , \label{state-linear}
\end{aligned}
\]
then it can be shown that the input-output mappings are additive:
\[
y_t = \psi_t(u_{0:t-1}, \delta_{0:t}) = \psi_t^u(u_{0:t-1}) + \psi_t^\delta(\delta_{0:t}) .
\]
In this case, the system is purifiable with
\[
\begin{aligned}
e_t &= \psi_t^\delta(\delta_{0:t}) \\
p_t(y_{0:t}, u_{0:t-1}) &= y_t - \psi_t^u(u_{0:t-1}) \\
q_t(e_{0:t}, u_{0:t-1}) &= e_t + \psi_t^u(u_{0:t-1}) .
\end{aligned}
\]

\section{Approximate solution methods}
\label{approximation}

Although Problem \eqref{OC-Q} is convex for convex systems, it remains challenging for two reasons. First, the decision variable $Q$ is infinite-dimensional. Second, the risk measures $\risk_0, \dots, \risk_J$ may be difficult to compute, or even ill-defined if some distributional information is lacking. We discuss methods for addressing the infinite-dimensionality in \S\ref{infinite-dimensionality} and for approximating risk measures in \S\ref{risk-approximation}. The upshot of this discussion is that several suboptimal control methods that perform well for linear systems can be applied to nonlinear convex systems using finite-dimensional convex optimization.

\subsection{Finite-dimensional restrictions}
\label{infinite-dimensionality}

\subsubsection{Open-loop model predictive control}
\label{OL-MPC}

A simple, effective method for overcoming the challenge of infinite dimensionality is open-loop MPC. In this method, at each time step we solve a version of the (finite-dimensional) open-loop optimal control problem \eqref{OC-OL} over a truncated, receding horizon. This generates a planned control input trajectory. We implement the first step in this plan, the system evolves, we update the state estimate and the process repeats. As discussed in \S\ref{convex-systems}, open-loop optimal control is convex for convex systems, even those with nonlinear dynamics. 

In one common variant of open-loop MPC, typically called certainty-equivalent MPC, the subproblem at each time step is solved under a single prediction of the disturbance trajectory. \cite{Mayne2000, Bemporad2002, Wang2010} This amounts to measuring risk with the \textit{predicted value} risk measure discussed in \S\ref{risk}. Other risk measures can be used in open-loop MPC, however, and can significantly improve performance. \cite{Schwarm1999, Blackmore2010, Cannon2011, Cinquemani2011, Schildbach2014, KircherSA2016, Campo1987, Magni2003, Mayne2005, CalafioreFagiano2013} Risk measures can be approximated if necessary, as discussed in \S\ref{risk-approximation}.

Open-loop MPC has several advantages. First, it often performs well in practice. Second, the design process is straightforward and intuitive; the primary design decisions are the prediction horizon, the terminal cost and/or constraints with which to augment the subproblems, and the algorithms for prediction and state estimation. A third advantage, and perhaps an underappreciated one, is that open-loop MPC admits very general uncertainty structures, including additive, multiplicative and others. In particular, since it does not require the Youla-type change of variables discussed in \S\ref{purifiable-Q}, open-loop MPC can be applied to systems that are not purifiable.

Open-loop MPC also has several disadvantages. First, its implementation is computationally intensive due to the use of online optimization; this can limit its scalability. (We note, however, that for linearly constrained linear systems with quadratic costs, open-loop MPC can be implemented efficiently using custom solvers that exploit the problem structure. \cite{Wang2010}) Second, open-loop MPC yields no closed-form expression for the control policy $\pi$; this complicates analysis of closed-loop stability, robustness, \etc \ Third, the method's performance can suffer somewhat due to the open-loop structure of the optimal control subproblems. This structure ignores the controller's opportunity to respond to future information as it becomes available, \ie, the controller's \textit{recourse}. The $Q$-design procedure partially addresses these disadvantages.

\subsubsection{Q-design}
\label{Q-design}

In open-loop MPC, the policy design procedure is straightforward (choosing a few parameters and subroutines), but implementation is computationally intensive due to the use of online optimization. In the $Q$-design procedure \cite{Kucera1980, Desoer1984, Anantharam1984, Boyd1991, Skaf2009, Skaf2010, Goulart2006}, by contrast, policy design is computationally intensive, but implementation is extremely efficient. In particular, no online optimization is needed. Another important distinction is that unlike MPC, $Q$-design yields a closed-form control policy. This facilitates analysis and simulation of the closed-loop system.

In $Q$-design, we design a suboptimal causal purified output feedback policy
\[
Q_{\theta} = \sum_{k=1}^K \theta_k Q^{(k)} .
\]
Here $Q^{(1)}, \dots, Q^{(K)} \in \mathcal Q$ are causal basis policies selected by the designer and $\theta = (\theta_1, \dots, \theta_K) \in \R^K$ is a parameter vector. The parameters are decided by solving a finite-dimensional analogue of Problem \eqref{OC-Q}, with $Q$ replaced by $Q_\theta$:
\bneq
\begin{array}{ll}
\underset{\theta}{\minimize}	& \risk_0 c_0(\phi(Q_{\theta}(\xi(\delta)), \delta), Q_{\theta}(\xi(\delta)), \delta)  \\ \label{OC-Qtheta} \tag{OC-$Q_\theta$}
\mbox{subject to}	& \risk_j c_j(\phi(Q_{\theta}(\xi(\delta)), \delta), Q_{\theta}(\xi(\delta)), \delta) \leq 0 , \\
&\qquad j = 1, \dots, J \\
				& Q_{\theta}(\xi(\delta)) = \sum_{k=1}^K \theta_k Q^{(k)}(\xi(\delta)) .
\end{array}
\eneq
Convex constraints on $\theta$ can also be added, \eg, to cultivate a particular structure in the control policy. Given $Q_\theta$, the unique corresponding causal output feedback policy can be constructed from the recursion in Theorem \ref{purifiability-Q-theorem}. 

Because $Q_\theta$ is affine in $\theta$ and convexity is preserved under composition with an affine mapping, the $Q$-design problem \eqref{OC-Qtheta} inherits the convexity of Problem \eqref{OC-Q} for convex systems. It follows that $Q$-design can be applied to convex systems, even nonlinear ones, using finite-dimensional convex optimization. Risk measures in Problem \eqref{OC-Qtheta} can be approximated if necessary, as discussed in \S\ref{risk-approximation}.

In principle, $Q$-design can solve the optimal control problem \eqref{OC-Q} to any degree of accuracy for a sufficiently large basis. In practice, $Q$-design has two disadvantages. First, its performance depends intimately on the choice of basis policies. In \cite{Skaf2010}, Skaf and Boyd explore the natural choice of affine policies in the context of linear systems. While affine $Q$-design performs well in a number of problems, the performance is typically not as good as open-loop MPC. To our knowledge, the general problem of finding good nonlinear basis policies has not yet been solved. A second disadvantage of $Q$-design is that, as discussed in \S\ref{purifiable-Q}, it can be applied only to purifiable systems. If the system is not purifiable, then an implementable output feedback policy may not be recoverable from the $Q$-designed policy.

\subsubsection{Closed-loop model predictive control}
\label{CL-MPC}

The final approximate solution method we discuss, closed-loop MPC, can be viewed as receding horizon $Q$-design. As in open-loop MPC, in closed-loop MPC we solve an optimal control subproblem at each time step over a truncated, receding horizon. The key difference is that the closed-loop MPC subproblems include a recourse model, \ie, a model of the controller's response to future information as it becomes available. More precisely, the open-loop MPC subproblems are truncated versions of the open-loop optimal control problem \eqref{OC-OL}, while the closed-loop MPC subproblems are truncated versions of the $Q$-design problem \eqref{OC-Qtheta}. Closed-loop MPC with affine recourse models is developed for linear systems in \cite{Skaf2010, Oldewurtel2008, Cinquemani2011, Schildbach2014, Goulart2007, CalafioreFagiano2013}.

Closed-loop MPC addresses the two disadvantages of $Q$-design. Like open-loop MPC, closed-loop MPC can be applied to systems that are not purifiable. To see this, we note that the first step of the recursion in Theorem \ref{purifiability-Q-theorem},
\[
\pi_0(y_0) = Q_0(y_0) ,
\]
can \textit{always} be implemented, even if later steps cannot. Closed-loop MPC also tends to be less sensitive than $Q$-design to the choice of basis policies; in closed-loop MPC the policy designed at each time step is merely a model of future recourse, while in $Q$-design the policy is the actual source of the implemented control inputs.

Closed-loop MPC can outperform open-loop MPC due to the inclusion of a recourse model. \cite{Skaf2010} Its design process is similarly straightforward; the only additional step is choosing basis policies. Like open-loop MPC, however, closed-loop MPC involves computationally intensive, online optimization. A closed-form expression for the closed-loop MPC policy is generally not available.

\subsection{Risk measure approximation}
\label{risk-approximation}

The control methods discussed in \S\ref{infinite-dimensionality} all require solving finite-dimensional convex optimization subproblems under uncertainty. These subproblems can be solved exactly in a few special cases. In general, however, we must resort to approximate solution methods. \cite{Prekopa1995, Shapiro2009, Ben-Tal2009} 

We now describe one approximate solution method based on sampling. This method can accommodate each of the risk measures discussed in \S\ref{risk}. We begin by obtaining samples $\delta^{(1)}, \dots, \delta^{(N)}$ from the distribution of the exogenous input trajectory $\delta$, or in the robust setting, from the corresponding uncertainty set. Samples could be obtained, \eg, from historical data or a pseudorandom number generator. Each risk measure $\risk_j$ is then replaced by an approximation $\hat \risk_j^N$ based on the samples $\delta^{(1)}, \dots, \delta^{(N)}$. This results in a convex optimization problem (with random input data) that can be solved using off-the-shelf software. Depending on the underlying risk measures, this sample-based approximate solution scheme is known as \textit{scenario optimization} \cite{Calafiore2005, Campi2008, Calafiore2010, Schildbach2013, Zhang2015, Campi2018} or \textit{sample-average approximation} \cite{Shapiro2009, Kleywegt2001, Homem-de-Mello2014, Kim2015}.

Table \ref{risk-approx-table} shows sample-based approximations for the five risk measures in Table \ref{risk-table}. The worst-case value and value-at-risk (which generate robust and chance constraints, respectively) can be approximated by maxima over the full sample. Expectation-type risk measures, including conditional value-at-risk, can be approximated using sample averages. Both of these approximations preserve the convexity of the underlying cost or constraint function. Theoretical bounds on violation probabilities for robust and chance constraints are available. \cite{Calafiore2005, Campi2008, Calafiore2010, Schildbach2013, Zhang2015, Campi2018} Some convergence results for sample-average approximation can be found in \cite{Shapiro2009, Kleywegt2001, Homem-de-Mello2014, Kim2015}. Variance reduction methods and hypothesis tests of solution quality can also be applied.  \cite{Shapiro2009}

\begin{table}
\centering
\caption{Approximations of the risk measures in Table \ref{risk-table} based on samples $X^{(1)}, \dots, X^{(N)}$ from the distribution of $X$.}
\smallskip
\begin{tabular}{l | l }
Risk measure, $\risk X$	&\begin{tabular}[x]{@{}l@{}} Sample-based approximation,\\ $\hat \risk^N(X^{(1)}, \dots, X^{(N)})$ \end{tabular} \\
\hline
Expected value			&$(X^{(1)} + \dots + X^{(N)})/N$ \\
Predicted value			&prediction based on $X^{(1)}, \dots, X^{(N)}$ \\
Worst-case value		&$\max \setof{X^{(1)}, \dots, X^{(N)}}$ \\
Value-at-risk			&$\max \setof{X^{(1)}, \dots, X^{(N)}}$ \\
\begin{tabular}[x]{@{}l@{}} Conditional \\\quad value-at-risk \end{tabular}	&\begin{tabular}[x]{@{}l@{}} $\inf_\alpha \{ \alpha + ( (X^{(1)}-\alpha)_+ + \dots$\\$\quad + \ (X^{(N)} - \alpha)_+)/(N(1-\beta)) \}$ \end{tabular}\\
\end{tabular}
\label{risk-approx-table}
\end{table}

To illustrate this method, we consider the open-loop optimal control problem \eqref{OC-OL}. The sample-based approximation to this problem is to
\[
\begin{array}{ll}
\underset{\tilde u}{\minimize}	& \hat \risk_0^N \big( c_0(\phi(\tilde u, \delta^{(1)}), \tilde u,  \delta^{(1)}), \dots, \\
&\quad c_0(\phi(\tilde u, \delta^{(N)}), \tilde u,  \delta^{(N)})  \big)\\ 
\mbox{subject to}			& \hat \risk_j^N \big( c_j(\phi(\tilde u, \delta^{(1)}), \tilde u,  \delta^{(1)}), \dots, \\
&\quad c_j(\phi(\tilde u, \delta^{(N)}), \tilde u,  \delta^{(N)})  \big) \leq 0 , \quad j = 1, \dots, J .
\end{array}
\]
The $Q$-design problem \eqref{OC-Qtheta} can be approximated similarly. The decision variables in these problems ($\tilde u$ and $\theta$, respectively) can be regularized to avoid overfitting the training data $\delta^{(1)}, \dots, \delta^{(N)}$. The computational complexity of the sampled problem generally scales linearly with $N$, so relatively large sample sizes can often be used. Different sample sizes $N_0, \dots, N_J$ can be used to approximate the different risk measures $\risk_0, \dots, \risk_J$.

\section{Numerical example}
\label{example}

In this section, we demonstrate the new capabilities developed in this paper through a numerical example. The example includes nonlinear dynamics, non-additive uncertainty and the $\cvar_\beta$ risk measure.

We consider a single-input, single-output system with dynamics
\bneq
x_{t+1} = (1 - w(\delta_{t+1})/20) (x_t)_+^2 - u_t + w(\delta_{t+1}) . \label{example-dynamics}
\eneq
We assume the controller has perfect state information ($y_t = x_t$), so that we can directly compare control methods in isolation from the state estimation problem. The unforced dynamics have a stable equilibrium at the origin and an unstable equilibrium at unity. Our primary control objective is to maintain the state within the basin of attraction of the stable equilibrium, \ie, to satisfy the constraint $c_1(x) \leq 0$, where
\[
c_1(x) = \max \setof{x_1, \dots, x_T} - 1 .
\]
We accept that this constraint may occasionally be violated due to the uncertain initial state and disturbances, but tolerate violations only if they are both infrequent and small. For this reason, we measure the risk of constraint violations by the conditional value-at-risk ($\risk_1 = \cvar_\beta$) at confidence level $\beta = 0.95$. This risk measure addresses both the probability and (conditionally) expected magnitude of constraint violations. As discussed in \S\ref{risk}, $\cvar_\beta c_1(x)$ can be constrained by minimizing a scalar $\alpha$ subject to
\[
\Expect (c_1(x) - \alpha)_+ \leq - \alpha (1-\beta) .
\]
We would also like to minimize the total cost of control effort,
\[
c_0(u,\delta) = p(\delta_1) \abs{u_0} + \dots + p(\delta_T) \abs{u_{T-1}} ,
\]
where $p(\delta_1), \dots, p(\delta_T) \geq 0$ are uncertain prices. We choose the cost risk measure $\risk_0 = \Expect$.

The dynamics \eqref{example-dynamics}, while nonlinear, are nondecreasing in $x_t$ and convex in $(x_t,u_t)$ for all $\delta_{t+1}$. (We assume $w(\cdot) \leq 20$, so the coefficient multiplying $(x_t)_+^2$ is nonnegative.) The cost function $c_0$ is convex in $u$ for all $\delta$, and the constraint function $c_1$ is convex and nondecreasing in each of $x_1, \dots, x_T$. By Theorem \ref{convex-system-theorem}, therefore, the system is convex. It is also purifiable, since the controller has perfect state information and the dynamics are invertible in the disturbance $w(\delta_{t+1})$. The purified output at time $t$ is $\xi_t(\delta_{0:t}) = w(\delta_t)$, with purifier
\[
p_t(y_{0:t},u_{0:t-1}) = \frac{ y_t - (y_{t-1})_+^2 + u_t }{ 1 + (y_{t-1})_+^2/20 } .
\]
The inverse purifier $q_t$ can be straightforwardly constructed from the system dynamics. Theorem \ref{purifiability-Q-theorem} therefore applies, establishing a bijection between disturbance feedback policies and state feedback policies for this system. This gives our (infinite-dimensional, convex) control design problem:
\bneq
\begin{array}{ll}
\underset{Q \in \mathcal Q, \ \alpha}{\minimize}	& \alpha + \Expect c_0(Q(\xi(\delta)),\delta)  \label{example-problem} \\
\mbox{subject to}			& \Expect (c_1(\phi(Q(\xi(\delta)),\delta)) - \alpha)_+ \leq - \alpha (1-\beta)  .
\end{array}
\eneq

We compare five suboptimal controllers for this system. The first is the optimal affine disturbance feedback policy, computed via the $Q$-design procedure. The second is a nonlinear $Q$-designed policy with piecewise quadratic basis policies tuned through trial and error. The third is the certainty-equivalent variant of open-loop MPC discussed in \S\ref{OL-MPC}; subproblems are solved under a single prediction of the exogenous input trajectory (in this case, its conditional expectation). The fourth is open-loop MPC with the same risk measures as Problem \eqref{example-problem}. The fifth is closed-loop MPC with the same risk measures and an affine disturbance feedback recourse model. We refer to the fourth and fifth controllers as open-loop and closed-loop \textit{scenario} MPC, respectively, as is common in the literature. \cite{Schildbach2014,Calafiore2013}

Implementing each of these controllers involves numerically solving subproblems generated by Problem \eqref{example-problem}. As discussed in \S\ref{risk-approximation}, we replace the expected values in these subproblems by sample averages over a training set of 1,000 sample exogenous input trajectories. We use a validation set of another 1,000 samples to tune parameters such as the basis policies in nonlinear $Q$-design and the MPC prediction horizon. We then compare the policies' performance in a test set of another 2,000 samples.

In simulations, the initial state has a half-normal distribution on $[1,\infty)$ with mean two. The disturbances $w(\delta_1), \dots, w(\delta_T)$ are identically Beta$(2,4)$ distributed, shifted and scaled to have mean zero and support $[-1,2]$. Each price $p(\delta_t)$ of control effort is exponentially distributed with mean $2+\cos(2\pi t/T)$; prices are lowest (on average) in the middle of the control horizon. The prices, disturbances and initial state are mutually independent.

Optimization is done in MATLAB using the Gurobi solver and the CVX modeling toolbox \cite{Grant2008}, which makes specifying our problems very easy. For example, the following code computes the optimal open-loop control input trajectory over $N$ sample exogenous input trajectories stored in \texttt{x0} ($1 \times N$), \texttt{w} ($T \times N$) and \texttt{p} ($T \times N$):
\begin{verbatim}
variables a u(T,1)
expression x(T,N)
x = stateTrajectories(x0,u,w,f,T);
minimize( a + sum(c0(u,p))/N )
subject to
sum(pos(c1(x) - a))/N <= -a*(1 - b)
\end{verbatim}
We wrote the vectorized function \texttt{stateTrajectories} to play the role of the input-state mapping $\phi$, recursively building the state trajectory in each of the $N$ Monte Carlo simulations from the dynamics function \texttt{f}. It typically speeds up CVX modeling by an order of magnitude compared to a loop over Monte Carlo runs. This function, along with all other code used in this paper, is online at \cite{KircherGithub}.

\begin{figure}
\centering
\includegraphics[width=\columnwidth]{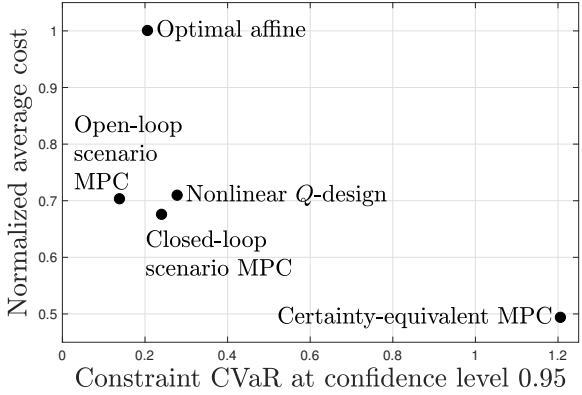}
\caption{Trade-offs between the cost of control effort and the risk of constraint violations in the test set of 2,000 sample exogenous input trajectories.}
\label{cost-robustness}
\end{figure}

Figure \ref{cost-robustness} shows the sample-average cost and constraint conditional value-at-risk of each controller in the test set of exogenous input trajectories. Table \ref{violationTable} shows the frequency with which each controller satisfies the constraint $x_t \leq 1$ in the test set. The five controllers strike different trade-offs between cost and robust constraint satisfaction. Certainty-equivalent MPC is aggressive, achieving low cost but poor robustness. The optimal affine controller satisfies constraints more robustly at the cost of increased control effort. The $Q$-designed nonlinear controller performs significantly better than the optimal affine controller, achieving similar cost and robustness to both of the scenario MPC variants. Including an affine recourse model in scenario MPC reduces cost slightly over open-loop scenario MPC, but also increases the risk of constraint violation. In this example, the robustness improvement of scenario MPC over certainty-equivalent MPC comes mainly from the use of samples, rather than the recourse model.

This numerical example demonstrates the key message of this paper: Suboptimal control methods that perform well for linear convex systems can be applied directly to nonlinear convex systems using the same software, and can perform similarly well. In this example, the $Q$-designed and predictive controllers show similar performance trends to the trends demonstrated for linear systems in \cite{Skaf2009, Skaf2010} and elsewhere. In particular, (a) scenario MPC is more robust than certainty-equivalent MPC, while achieving lower cost than affine controllers, and (b) nonlinear $Q$-designed controllers with good basis policies can be competitive with scenario MPC.

\section{Conclusion}
\label{conclusion}

\begin{table}
\centering 
\caption{Sample frequency of $x_t \leq 1$}
\begin{tabular}{l | l}
Controller	& Frequency \\
\hline
Optimal affine		& 99.5\% \\
\begin{tabular}[x]{@{}l@{}}Open-loop\\\quad scenario MPC\end{tabular}	& 99.3\% \\
Nonlinear $Q$-design		& 98.9\% \\
\begin{tabular}[x]{@{}l@{}}Closed-loop\\\quad scenario MPC\end{tabular}	& 98.7\% \\
\begin{tabular}[x]{@{}l@{}}Certainty-equivalent \\ \quad MPC\end{tabular}			& 74.6\% \\
\end{tabular}
\label{violationTable}
\end{table}

In this paper, we explored the class of nonlinear systems for which optimal control design under uncertainty can be cast as a convex optimization problem. We adopted a flexible approach to the risks associated with constraint violations and excessive costs, accommodating both the robust and stochastic views. We showed that open-loop optimal control is convex for nonlinear systems with convex or concave dynamics, provided the dynamics, cost and constraint functions satisfy certain monotonicity properties. We then showed that under the same conditions, closed-loop control design can be reformulated as a convex program if, in addition, the measured outputs can be reversibly `purified' of the influence of the past control inputs. 

The practical value of these results is the guarantee that, for a class of nonlinear systems, the subproblems solved in various suboptimal control methods are convex. This removes concerns about solvers getting stuck in local minima, and enables the use of convenient modeling software and reliable, efficient solvers. We illustrated this numerically for two methods, the $Q$-design procedure and model predictive control. 

There are a number of opportunities to extend this work. First, it would be interesting to survey applications in which nonlinear convex systems arise. We are aware of several examples, including cancer and HIV treatment scheduling \cite{Rantzer2014}, voltage control in electricity distribution networks \cite{Rantzer2014}, freeway traffic congestion management \cite{Schmitt2017} and energy storage control \cite{Khodabakhsh2016, Xu2017}, but this list is likely incomplete. 

Second, the requirement of purifiability appears to significantly restrict the uncertainty structures and nonlinearities admissible for $Q$-design. Exploring purifiability further could enable convex control design for a richer class of nonlinear systems. We note, however, that other methods, notably including model predictive control, do not require purifiability. 

Third, the class of nonlinear systems for which optimal control is convex could be further explored. We characterized only a subset of this class of systems; our method relied on the application of simple composition rules. While this method is constructive and compatible with convenient modeling software such as CVX \cite{Grant2008} and YALMIP \cite{Lofberg2004}, it is not exhaustive.

\appendices

\section{Proof of Theorem \ref{convex-system-theorem}}
\label{convex-system-proof}

Our task is to show that the function
\[
u \mapsto c_j((C(\delta)u + v(\delta),\eta(u,\delta)),u, \delta)
\]
is convex for all $\delta$. Here we have expressed the state trajectory as
\[
x = \bmat
x^{\rm aff} \\
x^{\rm cvx} \\
\emat = \bmat
C(\delta)u + v(\delta) \\
\eta(u,\delta) \\
\emat  .
\]
We will define the mappings $C$, $v$ and $\eta$ shortly.

We recall that convexity is preserved under composition with an affine mapping and under the composition of a convex nondecreasing function with a convex function \cite{Boyd2004}, and that $c_j$ is nondecreasing in each element of $x^{\rm cvx}$ by assumption. It therefore suffices to show that each row of $\eta(\cdot,\delta)$ is convex for all $\delta$.

To do this, we need the recursive definitions of $C$, $v$ and $\eta$. They are initialized by
\[
\begin{aligned}
C_1(\delta_{1}) &= B_1(\delta_1) \\
v_{1}(\delta_{0:1}) &= A_{1}(\delta_{1}) w_0(\delta_0) + w_{1}(\delta_{1}) \\
\eta_1(u_0,\delta_{0:1}) &= h_1(w_0(\delta_0), h_0(\delta_0), u_0, \delta_1) .
\end{aligned}
\]
For $t = 2, \dots, T$, 
\[
\begin{aligned}
C_{t}(\delta_{1:t}) &= \bmat A_{t}(\delta_{t}) C_{t-1}(\delta_{1:t-1}), & B_{t}(\delta_{t}) \emat \\
v_{t}(\delta_{0:t}) &= A_{t}(\delta_{t}) v_{t-1}(\delta_{0:t-1}) + w_{t}(\delta_{t}) \\
\eta_{t}(u_{0:t-1},\delta_{0:t}) &= h_{t}(C_{t-1}(\delta_{t-1}) u_{0:t-2} + v_{t-1}(\delta_{0:t-1}), \\
&\qquad \eta_{t-1}(u_{0:t-2},\delta_{0:t-1}), u_{t-1}, \delta_{t}) .
\end{aligned}
\]

The proof now proceeds by induction. At $t = 1$, $\eta_1$ is convex in $u_0$ by the convexity of $h_1$. For the inductive step, we suppose that each row of $\eta_{t-1}$ is convex in $u_{0:t-2}$. By assumption, each row of $h_t$ is nondecreasing in $x^{\rm cvx}_{t-1}$ and convex. Therefore, $\eta_t$ involves compositions of convex functions with the affine mapping $u_{0:t-2} \mapsto C_{t-1}(\delta_{t-1}) u_{0:t-2} + v_{t-1}(\delta_{0:t-1})$ and of convex nondecreasing functions with convex functions, namely the rows of $u_{0:t-2} \mapsto \eta_{t-1}(u_{0:t-2},\delta_{0:t-1})$. Applying our two composition rules, we see that each row of $\eta_t$ is convex. This concludes the inductive step and the proof.

\section{Proof of Theorem \ref{purifiability-Q-theorem}}
\label{purifiability-Q-proof}

Our task is to show that purifiability implies a bijection exists between causal $\pi$ in $y$ and causal $Q$ in some purified output that depends on $\delta$ alone. To do this, we suppose the system is purifiable. We build a causal $\pi$ in $y$ from a given causal $Q$ in the purified output
\[
e = \bmat
e_0 \\
e_1 \\
\vdots \\
e_{T-1} \\
\emat = \bmat
\xi_0(\delta_0) \\
\xi_1(\delta_{0:1}) \\
\vdots \\
\xi_{T-1}(\delta_{0:T-1}) \\
\emat =
\xi(\delta) .
\]
We recall that $e_0 = y_0$, always. The construction of $\pi$ from $Q$ is recursive:
\bneq
\begin{aligned}
\pi_0(y_0) &= Q_0(y_0) \\
\pi_t(y_{0:t}) &= Q_t(p_{0:t}(y_{0:t}, \pi_{0:t-1}(y_{0:t-1}))) , \\
&\qquad t = 1, \dots, T-1. \label{pi_recursion}
\end{aligned}
\eneq
Supposing this holds for $t = 0, \dots, \tau-1$, at $t = \tau$ we have
\[
\begin{aligned}
u_\tau &= Q_\tau(e_{0:\tau}) \\
&= Q_\tau(p_{0:\tau}(y_{0:\tau}, u_{0:\tau-1})) \\
&= Q_\tau(p_{0:\tau}(y_{0:\tau}, \pi_{0:\tau-1}(y_{0:\tau-1}))) \\
&\eqdef \pi_\tau(y_{0:\tau}) .
\end{aligned}
\]
By induction, therefore, the recursion \eqref{pi_recursion} holds for all $t = 0, \dots, T-1$. Given a causal $\pi$ in $y$, a causal $Q$ in $e$ can be constructed similarly from the recursion
\[
\begin{aligned}
Q_0(e_0) &= \pi_0(e_0) \\
Q_t(e_{0:t}) &= \pi_t(q_{0:t}(e_{0:t}, Q_{0:t-1}(e_{0:t-1})) ) , \\
&\qquad t = 1, \dots, T-1.
\end{aligned}
\]
This establishes the bijection and concludes the proof.

\section*{Acknowledgments}

The authors thank Anders Rantzer, Andy Ruina and Eilyan Bitar for their helpful comments, and gratefully acknowledge support from the National Science Foundation under grant 1711546.

\bibliography{bldg}

\begin{IEEEbiography}[{\includegraphics[width=1in,height=1.25in,clip,keepaspectratio]{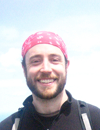}}]{Kevin J. Kircher} (M'15) is a mechanical engineering Ph.D. student at Cornell University. He received the B.S. degree in applied mathematics and physics with high honors from the University of Wisconsin-Milwaukee in 2008 and the M.Eng. degree in engineering physics from Cornell University in 2009. He then worked as a research associate at Lawrence Berkeley National Laboratory. His research interests include optimization and control under uncertainty, with applications to heating, cooling and power systems.
\end{IEEEbiography}

\begin{IEEEbiography}[{\includegraphics[width=1in,height=1.25in,clip,keepaspectratio]{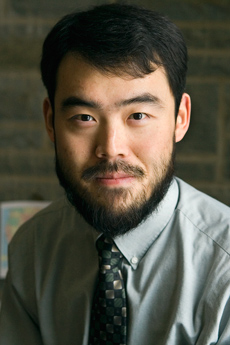}}]{K. Max Zhang} is an Associate Professor in the Sibley School of Mechanical and Aerospace Engineering at Cornell University, where he directs the Energy and Environment Research Laboratory. He received the Ph.D. degree in mechanical engineering from the University of California, Davis, in 2004. He was a visiting scientist at the United States Environmental Protection Agency (EPA) in 2000, 2002 and 2010--2012. Dr. Zhang received the EPA's Scientific and Technological Achievement Award and Cornell University's Engaged Scholar Prize.
\end{IEEEbiography}

\end{document}